\documentclass{elsart}

\usepackage{amsmath}
\usepackage{amssymb}
\usepackage[dvips]{graphicx}

\begin{document}

\begin{frontmatter}

\title{A backward Monte-Carlo method for solving parabolic partial
differential equations}

\author{Johan Carlsson\thanksref{email}}

\address{Oak Ridge National Laboratory,
P.O.~Box 2009, Oak Ridge, TN 37831--8071, USA}

\thanks[email]{E-mail: \textsf{carlssonja@ornl.gov}}

\begin{abstract}
A new Monte-Carlo method for solving linear parabolic partial differential
equations is presented.
Since, in this new scheme, the particles are followed \emph{backward} in
time, it provides great flexibility in choosing critical points in
phase-space at which to concentrate the launching of particles and thereby
minimizing the statistical noise of the sought solution.
The trajectory of a particle, $X_{i}(t)$, is given by the numerical
solution to the \emph{stochastic} differential equation naturally
associated with the parabolic equation.
The weight of a particle is given by the initial condition of the parabolic
equation at the point $X_{i}(0)$.
Another unique advantage of this new Monte-Carlo method is that it produces
a smooth solution, i.e.~without $\delta$-functions, by summing up the
weights according to the Feynman-Kac formula. \\[4mm]


\end{abstract}


\end{frontmatter}

\section{Introduction}

The Monte-Carlo method was conceived at the Los Alamos National Laboratory
during the Manhattan Project as a numerical method for solving the Boltzmann
equation governing the neutron distribution function in fissile material%
~\cite{mcm-first}. Since then it has found numerous other uses across many
fields of science. Overviews of its areas of applicability can be found in
any of a number of textbooks, see e.g.~Ref.~\cite{mcm-textbook}.
In plasma physics the Monte-Carlo method has a long history of being used for
solving the Fokker-Planck equation; see Ref.~\cite{mcm-transport} for some of
the earliest examples.
It has two main advantages: it keeps small the incremental effort of solving
a higher-dimensional problem, and it makes easy satisfying boundary conditions.
The Monte-Carlo method is typically worth consideration for three- or
higher-dimensional problems, or already in two dimensions if there are
internal boundary conditions imposed.
The biggest disadvantage is the unavoidable statistical noise caused by
the use of random numbers. This noise scales as the inverse square root
of the number of particles followed (and hence the number of arithmetic
operations and memory accesses required). The poor scaling is somewhat
offset by the fact that the Monte-Carlo method ``can be applied by many
computers working in parallel and independently'' as Metropolis and Ulam
pointed out more than half a century ago~\cite{mcm-first}.
In many cases, e.g.~when a tail distribution forms, low-density regions
of phase space are of particular interest. 
An example of such a case from plasma physics is when high-power radio waves
are launched into the plasma and absorbed through resonance with the gyration
of ions around magnetic field lines, resulting in the formation of a tail of
high-energy ions. The Monte-Carlo method has been used on
numerous occasions to solve the quasilinear Fokker-Planck equation, which
models such wave absorption~\cite{mcm-icrh}.
Due to the inverse-square-root scaling, the relative statistical error is
worst in exactly these interesting low-density regions.
Increasingly sophisticated weighting (i.e.~splitting of particles),
and reweighting, schemes~\cite{kasilov} have been suggested to make the
relative statistical error more constant throughout phase space. \\

The $\delta \! f$-method~\cite{collisionless-delta-f} also exists in a
collisional version%
~\cite{collisional-delta-f-1,collisional-delta-f-2a,collisional-delta-f-2b}.
There seem to be two different schools of thought on how to include the 
collisions: by making the particle trajectories
stochastic~\cite{collisional-delta-f-1} or by letting the collisions
enter the weight equation~%
\cite{collisional-delta-f-2a,collisional-delta-f-2b}.
With the latter approach the spreading of the weight causes a gradual
increase of the statistical error for a fixed number of particles. This makes
simulations problematic on timescales longer than a few collision times.
A potential cure for this particular ailment has recently been suggested%
~\cite{collisional-delta-f-2b}.
Since the $\delta \! f$-method assumes that the solution is a weak
perturbation of the equilibrium solution everywhere in phase space, it cannot
be used when the distribution function develops a tail. \\

The new Monte-Carlo method presented here is firmly based on the
well-established Feynman-Kac formula, which is briefly introduced
in section~\ref{sec:Feynman-Kac}. The Feynman-Kac formula puts the solution
of a parabolic partial differential equation on the form of a conditional
expectation value of a function of a stochastic variable.
\enlargethispage*{2.0ex}
In section~\ref{sec:mcm} it is shown how the numerical evaluation of this
expectation value takes the form of a Monte-Carlo method stepping backward
in time. Going backward in time allows us to choose the exact points in
phase space, e.g.~on an equidistant grid, at which we calculate the solution.
It also allows us to redistribute the statistical noise to regions of phase
space where it does the least harm. Finally, a comparison with the traditional
Monte-Carlo method can be found in section~\ref{sec:discussion}.

\section{The Feynman-Kac formula}

\label{sec:Feynman-Kac}

There is an intimate connection between linear parabolic partial
differential equations (PDEs) and stochastic differential equations
(SDEs). Let us illustrate by using the regular diffusion equation as
an example:
\begin{equation}
\frac{\partial f}{\partial t} =
\frac{\partial}{\partial x} \, D \, \frac{\partial f}{\partial x} \ ,
\quad x \in \mathbb{R} \ , \quad 0 \leq t \leq T \ , \label{eq:vanilla}
\end{equation}
with some initial condition $f(x,0) = \Phi(x)$.
Now, let the process $X(\tau)$ be governed by the SDE,
\begin{equation}
dX = \mu \, d\tau + \sigma \, dW \ , \label{eq:sde}
\end{equation}
where $W$ is a Wiener process, $\tau = T - t$, and we impose the
initial condition $X(0) = x$.
Using It\^o's formula [see Eq.~(\ref{eq:ito}) in appendix~\ref{app:ito}]
to differentiate $f(X(\tau),\tau)$, we get:
\begin{equation}
f(X(0),0) = f(X(T),T) + \int_{T}^{\, 0}
\left \{
\frac{\partial f}{\partial \tau} +
\mu \frac{\partial f}{\partial x} +
\frac{\sigma^{2}}{2} \frac{\partial^{2} \! f}{\partial x^{2}}
\right \}
d\tau +
\int_{T}^{\, 0} \sigma \, \frac{\partial f}{\partial x} \, dW \ .
\label{eq:proto}
\end{equation}
(To make this article more self-contained a brief introduction to stochastic
calculus and a mathematically somewhat questionable, but hopefully elucidating,
derivation of It\^o's formula is presented in appendix~\ref{app:ito} for those
who are unfamiliar with the formalism). By defining,
\begin{equation}
\mu = \frac{\partial D}{\partial x} \ , \quad
\sigma = \sqrt{2D} \ , \label{eq:musig}
\end{equation}
we make Eq.~(\ref{eq:sde}) the \emph{naturally associated} SDE
of the linear parabolic PDE~(\ref{eq:vanilla}).
Using the definition~(\ref{eq:musig}), Eq.~(\ref{eq:vanilla}), and the
initial conditions $f(x,t\!=\!0) = \Phi(x)$ and $X(\tau\!=\!0) = x$,
Eq.~(\ref{eq:proto}) becomes:
\begin{equation}
f(x,t\!=\!T) = \Phi(X(t\!=\!0)) +
\int_{T}^{\, 0} \sigma \, \frac{\partial f}{\partial x} \, dW \ . \nonumber
\end{equation}
Taking the expectation value of both sides we obtain:
\begin{equation}
f(x,t\!=\!T) = E[ \, \Phi(X(t\!=\!0)) \ | \ X(t\!=\!T) = x \, ] \ ,
\label{eq:feynman-kac}
\end{equation}
which is known as the stochastic Feynman-Kac representation of $f(x,T)$,
or the \emph{Feynman-Kac formula} for short. (See Ref.~\cite{ito} for a
thorough presentation of stochastic calculus and the Feynman-Kac formula).

\section{Backward Monte-Carlo method}

\label{sec:mcm}

In general, the expectation value on the RHS of Eq.~(\ref{eq:feynman-kac})
must be calculated numerically. We first integrate the SDE~(\ref{eq:sde}):
\begin{eqnarray}
\int_{\tau}^{\tau + \Delta t} dX &=&
\int_{\tau}^{\tau + \Delta t} \mu \, d\tau +
\int_{\tau}^{\tau + \Delta t} \sigma \, dW \ \Rightarrow \nonumber \\
X(\tau + \Delta t) - X(\tau) &=&
\mu [1 + \mathcal{O}(\sqrt{\Delta t})] \Delta t + \nonumber \\
&{}& \sigma [1 + \mathcal{O}(\sqrt{\Delta t})] [W(\tau + \Delta t) - W(\tau)]
\ .
\nonumber
\end{eqnarray}
The $\mathcal{O}(\sqrt{\Delta t})$ error terms come from the variation of
$\mu$ and $\sigma$ during the time step $\Delta t$. By the definition
of a Wiener process (see appendix~\ref{app:ito}), \linebreak
$W(\tau + \Delta t) - W(\tau) \in N(0,\sqrt{\Delta t})$ and we can write:
\begin{equation}
X(\tau + \Delta t) = X(\tau) + \mu \Delta t +
\zeta \sigma \sqrt{\Delta t} + \mathcal{O}(\Delta t) \ , \nonumber
\end{equation}
where $\zeta$ is a zero-mean, unit-variance Gaussian random number,
$\zeta \in N(0,1)$. The numerical approximation of the Feynman-Kac
formula~(\ref{eq:feynman-kac}) is simply:
\begin{equation}
%\boxed{
f(x,t\!=\!T) =  N^{-1} \sum_{i = 1}^{N} \Phi(X_{i}(t\!=\!0)) +
\mathcal{O}(\Delta t) + \mathcal{O}(N^{-1/2})%}
\ ,
\label{eq:numapprox}
\end{equation}
where $\mathcal{O}(N^{-1/2})$ is the statistical error and the stochastic
variables $X_{i}(t\!=\!0)$ are found by following the stochastic trajectories
given by:
\begin{equation}
%\boxed{
X_{i}(t - \Delta t) = X_{i}(t) + \mu \Delta t +
\zeta \sigma \sqrt{\Delta t} \ ,
\quad X_{i}(t = T) = x \, , \ i = 1 , \, \ldots , N %}
\ . \label{eq:motion}
\end{equation}

\begin{figure}[!t]
\begin{center}
\includegraphics[width=120mm]{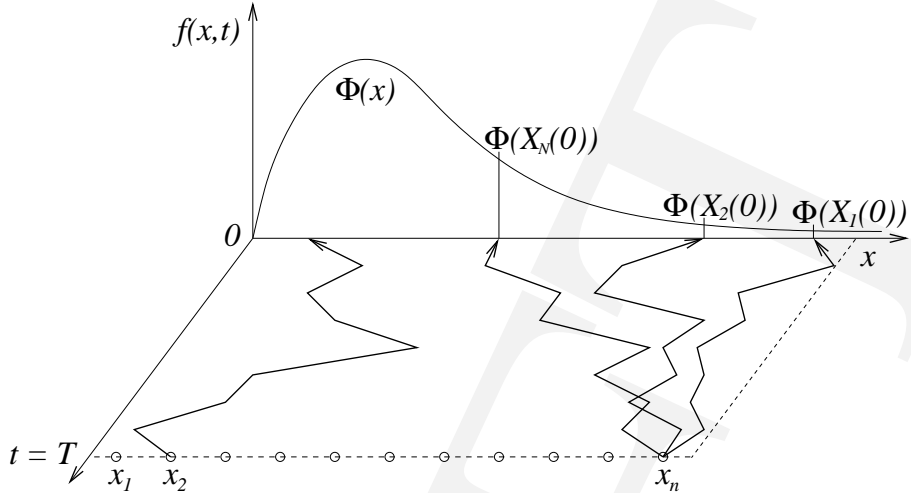}
\end{center}
\caption{The backward Monte-Carlo method. The trajectories as given by the
backward Monte-Carlo difference equation of motion~(\ref{eq:motion}).}
\label{fig:backward}
\end{figure}

A very simple algorithm, illustrated in Fig.~\ref{fig:backward}, can now be
used to solve Eq.~(\ref{eq:vanilla}).
Assume that we want the solution at time $t = T$ at the
points $x = x_{j} \, , \ j = 1 , \, \ldots , n$. Then, for each $j$, we simply
launch $N$ particles at $x = x_{j}$ and let them evolve according to the
backward Monte-Carlo equation of motion~(\ref{eq:motion}). As they reach
$t = 0$, we sample $\Phi (x)$ at their respective locations, $X_{i}(t\!=\!0)$,
and calculate the solution $f(x = x_{j}, t\!=\!T)$ as the average of the
sampled values, $\Phi (X_{i}(t\!=\!0))$, just as prescribed by
Eq.~(\ref{eq:numapprox}). \\

Notice that the relative statistical error should be roughly constant 
if $N$ is the same
for every $j$. Alternatively, as is indicated in Fig.~\ref{fig:backward},
we can concentrate the launching of particles to exactly those points where
a low-noise solution is desirable. In this sense, the backward Monte-Carlo
method offers a perfect weighting scheme. \\

The equivalents of Eqs.~(\ref{eq:numapprox}) and (\ref{eq:motion}) when
Eq.~(\ref{eq:vanilla}) is replaced by a more general linear parabolic PDE can
be found in appendix~\ref{app:general}.

\section{Discussion}

\label{sec:discussion}

The algorithm we arrived at in the previous section has a striking similarity
to the conventional Monte-Carlo method, but there are also some fundamental
differences.
When comparing the conventional Monte-Carlo difference equation of motion:
\begin{equation}
Y_{i}(t + \Delta t) = Y_{i}(t) + \mu \Delta t +
\zeta \sigma \sqrt{\Delta t} \ , \quad i = 1 , \, \ldots , N \ ,
\label{eq:motion-forward}
\end{equation}
to Eq.~(\ref{eq:motion}), it is evident that they both describe identical
trajectories, but for Eq.~(\ref{eq:motion}) these trajectories are traversed
backward in time. Since we are dealing with parabolic equations,
moving backward in time raises questions about time reversibility and the
change of entropy. The form of the solution~(\ref{eq:numapprox}) also
raises some suspicion; $f$ on the left-hand side is a macroscopic quantity and
so is $\Phi$ on the right-hand side, whereas $X_{i}(t\!=\!0)$ is microscopic.
In itself the backward Monte-Carlo difference equation of
motion~(\ref{eq:motion}) is perfectly legitimate; at the microscopic level
the motion is time-reversible since the entropy is undefined. The potential
danger lies in macroscopic information spilling over into the microscopic
world; i.e.~if the particles carried any information about the solution with
them going backward in time, then clearly the second law of thermodynamics
would be violated. Fortunately, the form of Eq.~(\ref{eq:numapprox})
guarantees that this will not happen since the particle weight
$\Phi(X_{i}(t\!=\!0))$ is undefined until $t = 0$.
Despite the superficial similarity between the Monte-Carlo difference
equations  of motion, (\ref{eq:motion}) and (\ref{eq:motion-forward}),
this is quite different from a conventional, forward Monte-Carlo method
where the particle weights are known at all times, $t \geq 0$. \\

\begin{figure}[!b]
\begin{center}
\includegraphics[width=120mm]{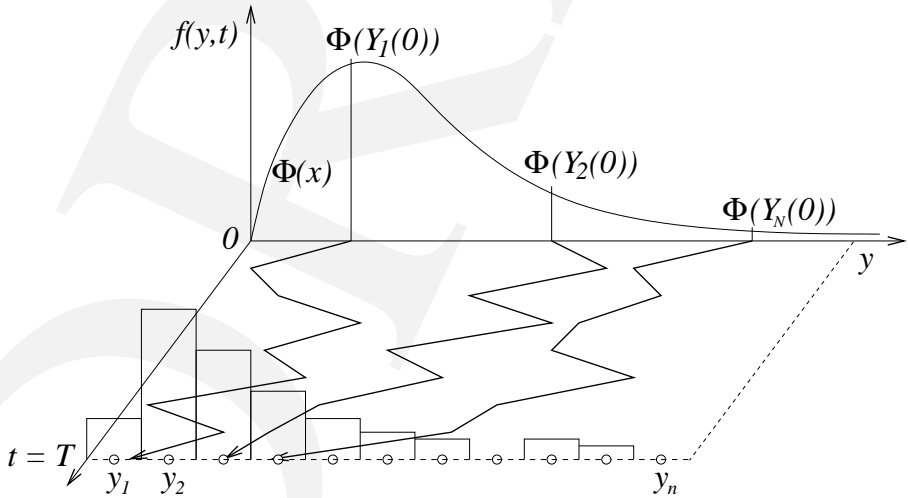}
\end{center}
\caption{Conventional, forward Monte-Carlo method with the same weights as the
backward method. The trajectories as given by the forward Monte-Carlo
difference equation of motion~(\ref{eq:motion-forward}).}
\label{fig:forward}
\end{figure}

Another difference between the backward and the forward Monte-Carlo
method is the very different character of the solutions. With the forward
Monte-Carlo method an obvious weighting scheme would be to use the same weights
as in the backward method and launch the particles with $Y_{i}(t\!=\!0)$
uniformly distributed over some sub-interval of $y$
(see Fig.~\ref{fig:forward}).
The solution then looks like
\begin{equation}
f(y,t\!=\!T) = N^{-1} \sum_{i=1}^{N} \Phi (Y_{i}(t\!=\!0)) \,
\delta (y - Y_{i}(t\!=\!T)) \ .
\nonumber
\end{equation}
Just to do something as simple as plotting the solution, the particles have
to be distributed in bins to smooth out the jaggedness and the solution
interpolated. This obviously means that a trade-off between resolution and
noisiness is unavoidable.
With the backward method, the solution [see Eq.~(\ref{eq:numapprox})]
is given as a numerical value at a point in phase space; it does not contain
any $\delta$-functions.
Points can be arbitrarily close together to give the desired resolution
without increasing the noise (see Fig.~\ref{fig:backward}).
It is also worth noting that the backward method makes it trivial to calculate
the solutions for a whole set of initial conditions,
$f(x,t\!=\!0) = \Phi_{m}(x)$, once the trajectories have been traced back to
$t = 0$ and the $X_{i}(t\!=\!0)$ are known. This makes the summing of samples
$\Phi(X_{i}(t\!=\!0))$ in Eq.~(\ref{eq:numapprox}) somewhat similar to the
convolution of a Green function with $\Phi(x)$. \pagebreak

In the forward Monte-Carlo method, the drift $\mu$ and the diffusion function
$\sigma$ are derived by taking moments of the single-particle distribution
function~\cite{moment}, and the forward stochastic Monte-Carlo difference
equation of motion~(\ref{eq:motion-forward}) is normally seen as something
rather artificial. From section~\ref{sec:mcm} we see that in the backward
Monte-Carlo method, the Monte-Carlo difference equation of
motion~(\ref{eq:motion}) has a more natural interpretation as the numerical
solution to the SDE~(\ref{eq:sde}) naturally associated with the parabolic
equation~(\ref{eq:vanilla}) that we wish to solve. \\

Finally, since it has a tendency to appear in contexts similar to this,
a very brief comment on the \emph{Langevin} equation, and why its use is
deprecated, is made in appendix~\ref{app:langevin}.

\section{Summary}

\label{sec:summary}

A backward Monte-Carlo method for solving parabolic differential equations
has been introduced. As compared to the conventional, forward Monte-Carlo
method, which is derived by taking moments of the single-particle distribution
function, the improved method originates from quite a different starting point:
the Feynman-Kac formula. \\

The stochastic Monte-Carlo difference equations of motion, including the
drift $\mu$ and the diffusion function $\sigma$, are identical in the
conventional and the new scheme, except for one vital difference: in the new
scheme the particles are followed backward in time. The similarity should
make it easy to retrofit the backward method to existing Monte-Carlo codes. \\

The solutions found with the forward Monte-Carlo method and the backward one,
however, take completely different forms. In the backward scheme, the solution
is smooth, unlike the jagged sum of $\delta$-functions associated with the
forward Monte-Carlo method. By default, the backward method also yields a
solution with a roughly constant relative statistical error throughout phase
space. In addition, it offers great flexibility in redistributing the
statistical noise to corners of phase space where it does minimal harm.
This latter capability makes the backward method particularly well suited
for cases where we are only interested in the solution in a small part of
phase space.

\begin{ack}

Prepared by the Oak Ridge National Laboratory, Oak Ridge, TN 37831--8071,
%managed by UT-BATTELLE for the U.S. DEPARTMENT OF ENERGY under contract
managed by UT--Battelle for the U.S.~Department of Energy under contract
DE--AC05--00OR22725.
Research was supported in part by an appointment to  the ORNL Postdoctoral
Research Associates Program, administered jointly by Oak Ridge National
Laboratory and the Oak Ridge Institute for Science and Education. \\

The author wishes to thank his colleagues in the Fusion Energy Division
\linebreak
Radio\-frequency Theory Group (Don Batchelor, Lee Berry, Mark Carter and
Fred Jaeger) for helpful comments during the work on this article.

\end{ack}

\appendix

\section{Stochastic calculus and It\^o's formula}

\label{app:ito}

We start with the stochastic differential equation~(\ref{eq:sde}),
repeated here for convenience:
\begin{equation}
dX = \mu \, d\tau + \sigma \, dW \ , \label{eq:sde3}
\end{equation}
where $W$ is a Wiener process. Ordinary stochastic
variables are just mappings from one probability space to another;
stochastic processes are time dependent. A stochastic process
$W(\tau), \tau \geq 0$ is a Wiener processes \emph{iff}:
\begin{itemize}
\item $W(0) = 0$
\item the increment $W(\tau+\Delta\tau) - W(\tau) \, , \ \Delta\tau > 0 \, ,$
is independent of $W(s) \, , \ s \leq \tau$
\item $W(\tau+\Delta\tau) - W(\tau) \in N \big (0,\sqrt{\Delta\tau} \, \big )$
\item $W(\tau)$ has continuous trajectories
\end{itemize}
We will first derive a differential identity that will be needed later
in this appendix. We first define $\Delta\tau = \tau_{j+1} - \tau_{j}$
and $\Delta W_{j} = W(\tau_{j+1}) - W(\tau_{j})$ with
$\tau_{j} = j \, \tau / n  \, , \ j = 0, 1, \ldots, n-1$. We are now ready
to introduce the stochastic variable $S_{n}(\tau)$,
\begin{equation}
S_{n}(\tau) = \sum_{j=0}^{n-1} ( \Delta W_{j} )^{2} \ . \label{eq:quad-var}
\end{equation}
If $dW/d\tau$ had existed, then clearly $S_{n}(\tau)$ would tend to zero as $n$
goes to infinity. But the derivative $dW/d\tau$ is undefined everywhere, so we
have to actually calculate the limit value. We take a congenially probabilistic
approach to this task. The expectation value $E[S_{n}]$ is trivial:
\begin{equation}
E[S_{n}] = \sum_{j=0}^{n-1} E \big [ ( \Delta W_{j} )^{2} \big ] =
\sum_{j=0}^{n-1} V [ \Delta W_{j} ] = \sum_{j=0}^{n-1} \Delta\tau_{j} = 
\sum_{j=0}^{n-1} ( \tau / n ) = \tau \ . \label{eq:mean}
\end{equation}
To establish that the expectation value (\ref{eq:mean}) is really the sought
limit of (\ref{eq:quad-var}), we must show that the variance $V[S_{n}]$ goes
to zero when $n$ goes to infinity.
We will start by calculating $E \big [ X^{4} \big ], X \in N(0,\sigma)$,
\begin{equation}
\begin{split}
E \big [ X^{4} \big ] = (2 \pi \sigma ^{2})^{-1/2} &
\int_{- \infty}^{\infty} x^{4} e^{ - \frac{ x^{2} }{ 2 \sigma ^{2} } } dx
= \cdots \\
= (2 \pi \sigma ^{2})^{-1/2} & \int_{- \infty}^{\infty}
3 \sigma ^{2} x^{2} e^{ - \frac{ x^{2} }{ 2 \sigma ^{2} } } dx =
3 \sigma ^{2} V[X] = 3 V[X]^{2} \ , \label{eq:help}
\end{split}
\end{equation}
where $(\cdots)$ is an integration by parts. With the help of the
identity~(\ref{eq:help}) we find
\begin{equation}
\begin{split}
V[S_{n}] =& \sum_{j=0}^{n-1} V \big [ ( \Delta W_{j} )^{2} \big ] =
\sum_{j=0}^{n-1} \Big ( E \big [ ( \Delta W_{j} )^{4} \big ] -
 E \big [ ( \Delta W_{j} )^{2} \big ]^{2} \Big ) \\
=& 2 \sum_{j=0}^{n-1} V [ ( \Delta W_{j} ]^{2} =
2 \sum_{j=0}^{n-1} ( \Delta\tau_{j} )^{2} = 
2 \sum_{j=0}^{n-1} ( \tau / n )^{2} = 2 \tau^{2} / n \ . \label{eq:variance}
\end{split}
\end{equation}
Now, since $V[S_{n}] \rightarrow 0 \, , \ n \rightarrow \infty$ and
$E[S_{n}] \rightarrow \tau \, , \ n \rightarrow \infty$, we will be bold
enough to draw the conclusion (inspired by the limit sum):
\begin{equation}
\lim_{n \rightarrow \infty} S_{n}(\tau) =
\lim_{n \rightarrow \infty} \sum_{j=0}^{n-1} ( \Delta W_{j} )^{2} = \tau \
\Rightarrow \ \int_{0}^{\tau} ( dW )^{2} =\int_{0}^{\tau} d\tau \ .
\label{eq:diffid}
\end{equation}
Now we are ready to calculate the differential $df(X(\tau),\tau)$ and start by
Taylor expanding $f$ to second order:
\begin{equation}
df = f_{x} \, dX + f_{\tau} \, d\tau + \half f_{xx} (dX)^{2} +
\half f_{\tau\tau} (d\tau)^{2} + f_{x\tau} \, dX \, d\tau \ . \nonumber
\end{equation}
Substituting Eq.~(\ref{eq:sde3}) for $dX$ and letting the
identity~(\ref{eq:diffid}) justify the ordering
$dW \gg d\tau = (dW)^{2} \gg dW d\tau \gg (d\tau)^{2}$, we get
\begin{equation}
df = \left \{ f_{\tau} + \mu f_{x} + \half \sigma^{2} f_{xx}
\right \} d\tau + \sigma f_{x} \, dW \ , \label{eq:ito}
\end{equation}
where only the two lowest orders ($dW$ and $d\tau$) have been kept. This is
the sought It\^o's formula, which is more rigorously  derived in
Ref.~\cite{ito}.

\newpage

\section{General linear parabolic PDE}

\label{app:general}

The backward Monte-Carlo method introduced in section~\ref{sec:mcm} can
be used to solve much more general linear parabolic PDEs than
Eq.~(\ref{eq:vanilla}). In this appendix we will generalize
Eqs.~(\ref{eq:numapprox}) and~(\ref{eq:motion}) to solve the following
equation:
\begin{equation}
\frac{\partial f}{\partial t} = \frac{\partial}{\partial x^{k}} \, D^{k\ell} \,
\frac{\partial f}{\partial x^{k}} +
\lambda f + S \ , \quad x^{k} \in \mathbb{R} \ , \quad 0 \leq t \leq T \ ,
\label{eq:general}
\end{equation}
with the initial condition $f(x^{k},0) = \Phi(x^{k})$, and
$D^{k\ell} = D^{k\ell}(x^{k}, t)$, $\lambda = \lambda(x^{k}, t)$, and
$S = S(x^{k}, t)$.
We will again try to find the Feynman-Kac representation of $f(x^{k},t)$, and
to do so we need the naturally associated SDEs. The matrix $D^{k\ell}$ is in
general not diagonal. In other words, the diffusion processes along the
different axes are in general correlated to some degree and the SDEs take the
form:
\begin{equation}
dX^{k} = \mu^{k} d\tau + A^{k\ell} dW^{\ell} \ . \nonumber
\end{equation}
It\^o's formula is trivial to generalize to multiple dimensions, and applying
it to $f$ we find the identities:
\begin{equation}
\mu^{k} = \frac{\partial D^{k\ell}}{\partial x^{\ell}}  \ , \nonumber
\end{equation}
and
\begin{equation}
A^{km} A^{\ell{}m} = 2 D^{k\ell} \ . \label{eq:algebraic}
\end{equation}
Finding the Feynman-Kac representation of the solution to 
Eq.~(\ref{eq:general}) is straightforward, with
$\Lambda(t) = \exp[\int_{0}^{t} \lambda(X^{k}(s), s) \, ds]$ we get
\begin{equation}
f(x^{k},T) = E \left [ \,
\Lambda(T) \, \Phi(X^{k}(0)) + \int_{0}^{T} \Lambda(t) \, S(X^{k}(t), t) \, dt
\ \Big | \ X^{k}(T) = x^{k} \, \right ] \ ,
\nonumber
\end{equation}
with the numerical approximation
\begin{equation}
f(x^{k},T) = N^{-1} \sum_{i = 1}^{N} \left \{
\Lambda(T) \, \Phi(X_{i}^{k}(0)) +
\int_{0}^{T} \Lambda(t) \, S(X_{i}^{k}(t), t) \, dt \right \} \ ,
\label{eq:general-numapprox}
\end{equation}
where
\begin{equation}
X_{i}^{k}(t - \Delta t) = X_{i}^{k}(t) + \mu^{k} \Delta t +
\zeta^{\ell} A^{k\ell} \sqrt{\Delta t} \ ,
\quad X_{i}^{k}(T) = x^{k} \, , \ i = 1 , \, \ldots , N \ .
\label{eq:general-motion}
\end{equation}
Here, $\zeta^{\ell}$ are uncorrelated, zero-mean, unit-variance Gaussian
random numbers, $\zeta^{\ell} \in N(0,1)$, and the matrix elements $A^{k\ell}$
solve the system of algebraic Eqs.~(\ref{eq:algebraic}). 

\section{Langevin equation}

\label{app:langevin}

A literature review on Monte-Carlo methods for solving parabolic equations
is impossible without occasionally coming across the Langevin
equation~\cite{langevin,chandrasekhar}:
\begin{equation}
\frac{d \mathbf{v}}{dt} = - \beta \mathbf{v} + \mathbf{A}(t) \ .
\label{eq:langevin}
\end{equation}
Here $\mathbf{v}$ is the velocity of a particle, and $\mathbf{A}(t)$ is
a ``fluctuating'' acceleration. The Langevin equation was historically used
to model Brownian motion~\cite{brownian}. \\

The following assumptions are being made about the ``fluctuating''
term $\mathbf{A}(t)$:
\begin{itemize}
\item $\mathbf{A}(t)$ is independent of $\mathbf{v}$.
\item $\mathbf{A}(t)$ varies extremely rapidly compared to the variations
of $\mathbf{v}$.
\end{itemize}
It should come as no surprise that the second assumption is problematic.
To quote Chandrasekhar~\cite{chandrasekhar}: ``But we should draw attention
even at this stage to the very drastic nature of assumptions implicit in the
very writing of an equation of the form~(\ref{eq:langevin}). For we have in
reality supposed that we can divide the phenomenon into two parts, one in
which the discontinuity of the events taking place is essential while in the
other it is trivial and can be ignored''. \\

Since the theory of stochastic calculus~\cite{ito} is on considerably firmer
mathematical footing, the SDE~(\ref{eq:sde}) should be allowed to supersede
the Langevin equation.

\end{document}